\documentclass{article}
\usepackage{amssymb}

\input{tcilatex}
\begin{document}

\textbf{Asymptotics of the number of partitions into p-cores and some
trigonometric sums.}

\[
Gert\text{ }Almkvist 
\]%
$\mid $

\textbf{1. Introduction.}

A p-core is a partition that has no hook numbers divisible by $p$. These
partitions occur in the representation theory of the symmetric group in
characteristic $p$.

Let $\ $\ $a_{p}(n)$ \ denote the number of partitions of \ $n$ \ into \
p-cores. Then we have the following generating function ( Garvan [7] , p.449)%
\[
\sum_{n=0}^{\infty }a_{p}(n)x^{n}=\dprod\limits_{j=1}^{\infty }\frac{%
(1-x^{pj})^{p}}{(1-x^{j})}
\]%
Garvan also gives the following asymptotic formula for \ $a_{p}(n)$ \ if $\ p
$ \ is a prime $\geq 5$%
\[
a_{p}(n-\frac{p^{2}-1}{24})=K_{p}\sum_{d\mid n}\QOVERD( ) {d}{p}(\frac{n}{d}%
)^{(p-3)/2}+O(n^{(p-3)/4+\varepsilon })
\]%
where \ $\QOVERD( ) {d}{p}$ \ is the Legendre symbol and \ $K_{p}$ is a
constant depending only on \ $p$. \ Implicitly in Garvan's paper one has \ $%
K_{p}=1/c_{p}$ \ where \ $c_{p}$ \ seems to be an integer. This is proved
here using results from [2] and [7]

Using the rather obscure methods described in [1] we will find another
asymptotic formula for \ $a_{p}(n)$ . Deleting some, usually small terms, it
is proved that we get the same formula as Garvan. For this we have to use
that%
\[
\frac{\eta (\tau )^{p}}{\eta (p\tau )} 
\]%
is a modular form for the group $\Gamma _{0}(p)$ \ with character \ $%
\QOVERD( ) {d}{p}$ \ (Ogg [9], p.28).

As a byproduct we obtain some trigonometric sums which turn out to be
integer valued. Here is a sample%
\[
\frac{\sqrt{p}}{2^{r+1}}\sum_{j=1}^{(p-1)/2}\cot ^{(r)}(\frac{j^{2}\pi }{p}%
)=-(-1)^{[(r+1)/2]}\QOVERD( ) {-1}{p}\frac{p^{r+1}}{r+1}%
\sum_{j=1}^{(p-1)/2}B_{r+1}(frac(\frac{j^{2}}{p})) 
\]%
are integers. Here \ $B_{n}(x)$ \ is the Bernoulli polynomial and \ $%
frac(x)=x-[x]$ \ is the fractional part of \ $x.$\ 
\[
\]

\textbf{2. An asymtotic formula for \ }$\mathbf{a}_{\mathbf{p}}\mathbf{(n).}$

Let 
\[
F(x)=\dprod\limits_{j=1}^{\infty }(1-x^{j})^{-1}
\]%
be the generating function for \ $p(n)$, the number of partitions of \ $n.$
Then%
\[
f(x)=\dprod\limits_{j=1}^{\infty }\frac{(1-x^{pj})^{p}}{1-x^{j}}=\frac{F(x)}{%
F(x^{p})^{p}}
\]%
We want to study \ $f(x)$ \ near the point%
\[
x=\exp (2\pi ih/k)
\]%
where \ $(h,k)=1.$ But it is well kown that \ (see Apostol [5], p.106)%
\[
F(\exp (2\pi ih/k-t)
\]%
\[
=\exp (\pi is(h,k))(\frac{kt}{2\pi })^{1/2}\exp (\frac{\pi ^{2}}{6k^{2}}%
t^{-1}-\frac{t}{24})F(\exp (\frac{2\pi iH}{k}-\frac{4\pi ^{2}}{k^{2}t})
\]%
where%
\[
Hh\equiv -1\text{ }\func{mod}\text{ }k
\]%
and%
\[
s(h,k)=\sum_{j=1}^{k-1}((j/k))((jh/k))
\]%
is the classical Dedekind sum. Here \ $((x))=x-[x]-1/2$ \ if \ $x\notin 
\mathbf{Z}$ \ and \ $=0$ \ if \ $x\in \mathbf{Z}$ . It follows that%
\[
f(\exp (2\pi ih/k-t))=\exp (\pi i(s(h,k)-ps(ph,k))
\]%
(*)%
\[
\cdot (\frac{2\pi }{k})^{(p-1)/2}p^{-p/2}t^{-(p-1)/2}\exp (\frac{p^{2}-1}{24}%
t)\frac{F(\exp (2\pi iH/k-4\pi ^{2}/k^{2}t))}{F(\exp (2\pi iB/k-4\pi
^{2}/k^{2}pt))^{p}}
\]%
where%
\[
Bph\equiv -1\text{ }\func{mod}\text{ }k\text{ \ if \ }(k,p)=0
\]%
If \ $(k,p)\neq 1$ \ we put \ $B=0.$ Define%
\[
A_{p}(k,n)=\sum_{(h,k)=1}\exp (\pi i(s(h,k)-ps(ph,k))-2\pi ihn/k)
\]%
Neglecting all analytic difficulties we compute \ $a_{p}(n)$ \ as the n-th
Fourier coefficient of \ $f(\exp (iy))$%
\[
a_{p}(n)=\frac{1}{2\pi }\dint\limits_{-\pi }^{\pi }f(\exp (iy))\exp (-iny)dy
\]%
Put \ $y=2\pi h/k+\varphi $ \ and sum over all \ $h=1,2,...k-1$ \ with \ $%
(h,k)=1$ \ (we also replace the interval \ $(-\pi ,\pi )$ \ with \ $(-\infty
,\infty )$ \ in order to get a Fourier transform). Let us also asume that \
(*) \ is valid when we put \ $t=-i\varphi $ . We replace the last factor in
(*) by \ $1.$ Then the contribution from the singular points \ $\exp (2\pi
ih/k)$ \ with \ $(h,k)=1$ \ will be%
\[
\Phi _{k}(n)\approx (\frac{2\pi }{k})^{(p-1)/2}p^{-p/2}A_{p}(k,n)\frac{1}{%
2\pi }\dint\limits_{-\infty }^{\infty }(-i\varphi )^{-(p-1)/2}\exp
(-(n+(p^{2}-1)/24)i\varphi )d\varphi 
\]%
\[
=(\frac{2\pi }{k})^{(p-1)/2}p^{-p/2}A_{p}(k,n)\frac{%
(n+(p^{2}-1)/24)^{(p-3)/2}}{((p-3)/2)!}
\]%
(the integral exists in the distributional sense). Summing over \ $k$ \ we
obtain%
\[
a_{p}(n)\approx (\frac{2\pi }{k})^{(p-1)/2}p^{-p/2}A_{p}(k,n)\frac{%
(n+(p^{2}-1)/24)^{(p-3)/2}}{((p-3)/2)!}\sum_{k=1}^{\infty }\frac{A_{p}(k)}{%
k^{(p-1)/2}}
\]%
which looks rather different from Garvan's formula (to begin with, it is an
infinite sum!).

Making some numerical computations we found the following

\textbf{1. }The formula gives about half of the digits correctly. E.g. 
\[
a_{17}(1000)=1829\text{ }06764\text{ }82504 
\]%
while the approximation gives%
\[
1829\text{ }06768\text{ }71721 
\]

\textbf{2. }The \ $A_{p}(k,n)$ \ are integers if \ $(k,p)=1.$ Computing many
values of \ $A_{p}(k,n)$ \ I finally came up with%
\[
\]

\textbf{Conjecture 1. \ }If \ $p\geq 5$ \ is prime and \ $(k,p)=1$ \ then%
\[
A_{p}(k,n)=\QOVERD( ) {k}{p}c_{k}(n+\frac{p^{2}-1}{24}) 
\]%
where%
\[
c_{k}(n)=\sum_{(h,k)=1}\exp (2\pi ih/k) 
\]%
is the Ramanujan sum (known to be an integer). This, however, follows from
the more precise 
\[
\]

\textbf{Conjecture 2. }If \ $p\geq 5$ \ is prime and \ $(k,p)=1$ \ then%
\[
ps(ph,k)-s(h,k)=\frac{(p^{2}-1)h}{12k}+\func{integer} 
\]%
where the integer is even if and only if \ $\QOVERD( ) {k}{p}=1$%
\[
\]

We will prove Conjecture 2 later. Let us first see what we can get from our
asymptotic formula from Conjecture 1. Since we know nothing about \ $%
A_{p}(k,n)$ \ when \ $k\equiv 0$ \ $\func{mod}$ $p$ , we throw away those
terms. Then we obtain%
\[
a_{p}(n-\frac{p^{2}-1}{24})\approx \frac{(2\pi )^{(p-1)/2}p^{-p/2}}{%
((p-3)/2)!}\sum_{k=1,p\nmid k}A_{p}(k-\frac{p^{2}-1}{24})\frac{n^{(p-3)/2}}{%
k^{(p-1)/2}} 
\]%
\[
=\frac{(2\pi )^{(p-1)/2}p^{-p/2}}{((p-3)/2)!}\sum_{k=1,k\nmid p}\QOVERD( )
{k}{p}\frac{c_{k}(n)}{k^{(p-1)/2}} 
\]%
Ramanujan ([10]) computed the sum%
\[
\sum_{k=1}^{\infty }\frac{c_{k}(n)}{n^{1+s}} 
\]%
but the factor \ $\QOVERD( ) {k}{p}$ \ causes some difficulties.

\textbf{Proposition 2.1. }\ We have the formula%
\[
\sum_{k=1}^{\infty }\QOVERD( ) {k}{p}\frac{c_{k}(n)}{n^{1+s}}=\frac{%
p^{1+s}\sum_{d\mid n}\QOVERD( ) {d}{p}d^{-s}}{\sum_{j=1}^{p-1}\QOVERD( )
{j}{p}\zeta (1+s,j/p)} 
\]%
where%
\[
\zeta (s,a)=\sum_{n=0}^{\infty }\frac{1}{(n+a)^{s}} 
\]%
is Hurwitz's $\zeta -$function.

\textbf{Proof:} First we note that 
\[
c_{k}(n)=\sum_{d\mid (n,k)}d\mu (k/d) 
\]%
where $\ \mu $ \ is the M\"{o}bius function. It follows%
\[
\sum_{k=1}^{\infty }\QOVERD( ) {k}{p}\frac{c_{k}(n)}{n^{1+s}}%
=\sum_{k=1}^{\infty }\QOVERD( ) {k}{p}\frac{1}{n^{1+s}}\sum_{d\mid
(n,k)}d\mu (k/d) 
\]%
\[
=\sum_{d\mid n}\sum_{j=1}^{\infty }\QOVERD( ) {jd}{p}\frac{1}{(jd)^{1+s}}%
d\mu (j)=\sum_{d\mid n}\QOVERD( ) {d}{p}\frac{1}{d^{s}}\sum_{j=1}^{\infty
}\QOVERD( ) {j}{p}\frac{\mu (j)}{j^{1+s}} 
\]%
But \ $\QOVERD( ) {j}{p}$ \ and \ $\frac{1}{j^{1+s}}$ \ are both completely
multiplicative functions of \ $j$ \ and we get%
\[
\sum_{j=1}^{\infty }\mu (j)\QOVERD( ) {j}{p}\frac{1}{j^{1+s}}=\frac{1}{%
\sum_{j=1}^{\infty }\QOVERD( ) {j}{p}\frac{1}{j^{1+s}}}=\frac{p^{1+s}}{%
\sum_{j=1}^{p-1}\QOVERD( ) {j}{p}\zeta (1+s,j/p)} 
\]%
and this completes the proof.%
\[
\]

Returning to the asymptotic formula we obtain%
\[
a_{p}(n)\approx K_{p}\sum_{d\mid n}\QOVERD( ) {d}{p}(n/d)^{(p-3)/2} 
\]%
which agrees with Garvans formula with%
\[
K_{p}=\frac{(2\pi )^{(p-1)/2}}{\sqrt{p}((p-3)/2)!\sum_{j=1}^{p-1}\QOVERD( )
{j}{p}\zeta ((p-1)/2,j/p)} 
\]%
\[
\]

\textbf{Proposition 2.2. }We have \ $K_{p}=1/c_{p}$ \ where

(i) \ \ 
\[
c_{p}=\frac{\sqrt{p}((p-3)/2)!}{(2\pi )^{(p-1)/2}}\sum_{j=1}^{p-1}\QOVERD( )
{j}{p}\zeta ((p-1)/2,j/p) 
\]

(ii)%
\[
c_{p}=\frac{1}{2}\QOVERD( ) {-2}{p}p^{(p-1)/2}\sum_{j=1}^{p-1}\QOVERD( )
{j}{p}\zeta (-(p-3)/2,j/p) 
\]

(iii)%
\[
c_{p}=(-1)^{(p-3)/2}\frac{\sqrt{p}}{2^{(p+1)/2}}\sum_{j=1}^{p-1}\QOVERD( )
{j}{p}\cot ^{((p-3)/2)}(\pi j/p) 
\]

(iv)%
\[
c_{p}=-\QOVERD( ) {-2}{p}\frac{p^{(p-1)/2}}{p-1}\sum_{j=1}^{p-1}\QOVERD( )
{j}{p}B_{(p-1)/2}(j/p) 
\]%
is a positive integer ( $p\geq 5$ \ prime ). \ If \ \ $p\equiv 3$ mod $4$ \
we also have the formulas

(v)%
\[
c_{p}=\frac{\sqrt{p}}{2^{(p+1)/2}}\sum_{j=1}^{p-1}\cot ^{((p-3)/2)}(\frac{%
\pi j^{2}}{p}) 
\]

(vi)%
\[
c_{p}=-\QOVERD( ) {-2}{p}\frac{p^{(p-1)/2}}{p-1}%
\sum_{j=1}^{p-1}B_{(p-1)/2}(frac(j^{2}/p))
\]%
\[
\]%
\[
\text{\textbf{Table of }}\mathbf{c}_{\mathbf{p}}
\]%
\[
\begin{tabular}{|l|l|}
\hline
$p$ & $c_{p}$ \\ \hline
5 & 1 \\ \hline
7 & 8 \\ \hline
11 & 1275 \\ \hline
13 & 33463 \\ \hline
17 & 5999 01794 \\ \hline
19 & 37084 43635 \\ \hline
23 & 2753 39898 05352 \\ \hline
29 & 66758 49413 21255 71317 \\ \hline
31 & 121 29134 29668 98388 66288 \\ \hline
\end{tabular}%
\]%
\[
\]

\textbf{Proof: }These formulas are conseqences of more general identities.
We will use the Finite Fourier Transformation \ (FFT) ( $k>1$ a fixed
integer )%
\[
\widehat{f}(\mu /k)=\sum_{j=0}^{k-1}f(j/k)\exp (-2\pi ij\mu /k) 
\]%
We make a little table of FFT for our purposes%
\[
\begin{tabular}{|l|l|}
\hline
$f$ & $\widehat{f}$ \\ \hline
$\zeta (s,x)$ & $k^{s}l(s,1-x)$ \\ \hline
$\QOVERD( ) {x}{k}$ & $(-i)^{((k-1)/2)^{2}}\sqrt{k}\QOVERD( ) {x}{k}$ if $\
k $ \ is prime \\ \hline
$B_{n}(x)$ & $kn(-2ki)^{-n}\cot ^{(n-1)}(\pi x)$ \\ \hline
\end{tabular}%
\]%
where%
\[
l(s,x)=\sum_{n=1}^{\infty }\frac{\exp (2\pi ix)}{n^{s}} 
\]%
is the periodic $\zeta -$function. There seems to be no reference to the
first and third formula. Therefore they are proved in an Appendix. The
second formula is in Apostol [4] p.195.

Define a scalar product%
\[
\left\langle f,g\right\rangle =\sum_{j=0}^{k-1}f(j/k)\overline{g(j/k)} 
\]%
Then we have Parseval's formula%
\[
\left\langle \widehat{f},\widehat{g}\right\rangle =k\left\langle
f,g\right\rangle 
\]%
Taking 
\[
f(j/p)=\QOVERD( ) {j}{p} 
\]%
and%
\[
g(j/p)=\zeta s,j/p) 
\]%
we get%
\[
\sum_{j=1}^{p-1}\QOVERD( ) {j}{p}\zeta (s,j/p)=\frac{1}{p}%
(-i)^{((p-1)/2)^{2}}\sqrt{p}p^{s}\sum_{j=1}^{p-1}\QOVERD( ) {j}{p}\overline{%
l(s,1-j/p)} 
\]%
In particular, let \ $s=r+1$ \ be a positive integer:%
\[
U=\sum_{j=1}^{p-1}\QOVERD( ) {j}{p}\zeta (r+1,j/p)=(-i)^{((p-1)/2)^{2}}\sqrt{%
p}p^{r}\sum_{j=1}^{p-1}\QOVERD( ) {j}{p}\sum_{n=1}^{\infty }\frac{\exp (2\pi
inj/p)}{n^{r+1}} 
\]

\textbf{Case 1: }$p\equiv 1$ $\func{mod}$ $4$

Then \ $\QOVERD( ) {-1}{p}=1$ \ and \ $(-i)^{((p-1)/2)^{2}}=1.$ Hence%
\[
U=\sqrt{p}p^{r}\sum_{j=1}^{p-1}\QOVERD( ) {j}{p}\sum_{n=1}^{\infty }\frac{%
\cos (2\pi nj/p)}{n^{r+1}} 
\]%
\[
=\frac{\sqrt{p}p^{r}(2\pi )^{r+1}(-1)^{(r+3)/2}}{2(r+1)!}\sum_{j=1}^{p-1}%
\QOVERD( ) {j}{p}B_{r+1}(j/p) 
\]%
if r is odd (see Apostol [3], p.267). But $\ r=(p-3)/2$ \ is odd so we get 
\[
K_{p}=-\QOVERD( ) {2}{p}\frac{p-1}{p^{(p-1)/2}\sum_{j=1}^{p-1}\QOVERD( )
{j}{p}B_{(p-1)/2}(j/p)} 
\]%
\[
\]

\textbf{Case 2: }$p\equiv 3$ $\func{mod}$ $4$

Then \ $\QOVERD( ) {-1}{p}=-1$ \ so \ $\QOVERD( ) {p-j}{p}=-\QOVERD( )
{j}{p} $ \ and \ $(-i)^{((p-1)/2)^{2}}=-i$ . Hence%
\[
U=-i\sqrt{p}p^{r}\sum_{j=1}^{p-1}\QOVERD( ) {j}{p}\sum_{n=1}^{\infty }\frac{%
i\sin (2\pi nj/p)}{n^{r+1}} 
\]%
\[
=\frac{\sqrt{p}p^{r}(2\pi )^{r+1}(-1)^{(r+2)/2}}{2(r+1)!}\sum_{j=1}^{p-1}%
\QOVERD( ) {j}{p}B_{r+1}(j/p) 
\]%
if r is even. In particular \ $r=(p-3)/2$ \ we get%
\[
K_{p}=-\QOVERD( ) {-2}{p}\frac{p-1}{p^{(p-1)/2}\sum_{j=1}^{p-1}\QOVERD( )
{j}{p}B_{(p-1)/2}(j/p)} 
\]%
(one has to check that the sign is correct). The rest of the proof will
follow from the results in the next section.%
\[
\]

\textbf{3. Some trigonometric sums.}

\textbf{Proposition 3.1. }Define for a prime \ $p\geq 5$ \ 
\[
S(r,p)=-\QOVERD( ) {-1}{p}\frac{\sqrt{p}}{2^{r+2}}\sum_{j=1}^{p-1}\QOVERD( )
{j}{p}\cot ^{(r)}(\pi j/p) 
\]%
and%
\[
T(r,p)=(-1)^{[(r-1)/2]}\frac{p^{r+1}}{2(r+1)}\sum_{j=1}^{p-1}\QOVERD( )
{j}{p}B_{r+1}(j/p) 
\]%
Then for \ $(p,r+1)=1$ \ both \ $S(r,p)$ \ and \ $T(r,p)$ \ are nonnegative
integers and equal.

\textbf{Proof: }By symmetry one observes that \ $S(r,p)=T(r,p)=0$ \ unless%
\[
p\equiv 1\text{ }\func{mod}\text{ }4\text{ \ and \ }r\text{ \ odd} 
\]%
or%
\[
p\equiv 3\text{ }\func{mod}\text{ }4\text{ \ and \ }r\text{ \ even} 
\]%
Using Parseval's formula we obtain 
\[
T(r,p)=(-1)^{[(r-1)/2]}\frac{p^{r+1}}{2(r+1)}\sum_{j=1}^{p-1}\QOVERD( )
{j}{p}B_{r+1}(j/p) 
\]%
\[
=(-1)^{[(r-1)/2]}\frac{p^{r+1}}{2(r+1)}\sum_{j=1}^{p-1}(-i)^{((p-1)/2)^{2}}%
\sqrt{p}\QOVERD( ) {j}{p}\overline{\frac{p(r+1)}{(-2pi)^{r+1}}\cot
^{(r)}(\pi j/p)} 
\]%
\[
=-\QOVERD( ) {-1}{p}\frac{\sqrt{p}}{2^{r+2}}\sum_{j=1}^{p-1}\QOVERD( )
{j}{p}\cot ^{(r)}(\pi j/p)=S(r,p) 
\]%
The sign is easily checked in the two cases mentioned in the beginning of
the proof.

That \ $T(r,p)$ \ is an integer follows from the following facts

(i) \ $k^{n}(B_{n}(j/k)-B_{n})$ \ is an integer (Almkvist-Meurman [2] ).
Note that \ $\sum_{j=1}^{p-1}\QOVERD( ) {j}{p}B_{n}=0$

(ii) \ If \ $(k,n)=1$ \ then this integer is divisible by \ $k$ \
(Dokshitzer [6], p.3)

(iii)%
\[
\sum_{j=1}^{p-1}\QOVERD( ) {j}{p}B_{r+1}(j/p)=2\sum_{j=1}^{(p-1)/2}\QOVERD(
) {j}{p}B_{r+1}(j/p) 
\]%
in both cases.

It remains to show that, say \ $S(r,p)$ \ is nonnegative. The idea of proof
is to show that the first term \ $\left\vert \cot ^{(r)}(\pi /p)\right\vert $
dominates the sum.

\textbf{(a) }$\left\vert \cot ^{(r)}(x)\right\vert $ \ is decreasing in the
interval \ $(0,\pi /2)$

This follows from that \ there is a polynomial \ $f_{r}(t)$ \ of degree \ $%
r+1$ \ such that \ $\left\vert \cot ^{(r)}(x)\right\vert =f_{r}(\cot (x))$ \
\ and all the coefficints of \ $f_{r}$ \ are nonnegative. Indeed%
\[
f_{1}(t)=1+t^{2} 
\]%
and 
\[
f_{r+1}(t)=(1+t^{2})f_{r}^{\prime }(t) 
\]

\textbf{(b) \ }$r=1.$ Then \ $\left\vert \cot ^{(1)}(x)\right\vert =\frac{1}{%
\sin ^{2}(x)}$ \ and we have%
\[
\frac{1}{\sin ^{2}(\pi /p)}\geq \frac{p^{2}}{\pi ^{2}} 
\]%
We will show that 
\[
\sum_{j=2}^{(p-1)/2}\frac{1}{\sin ^{2}(\pi j/p)}<\frac{p^{2}}{\pi ^{2}} 
\]%
if \ $p$ \ is large enough. One has the inequality (see [3], p.45)%
\[
\frac{1}{\sin (x)}\leq \frac{1}{x}+\frac{x}{3}\text{ \ if \ }0<x\leq \frac{%
\pi }{2} 
\]%
Hence%
\[
\sum_{j=2}^{(p-1)/2}\frac{1}{\sin ^{2}(\pi j/p)}\leq \frac{p^{2}}{\pi ^{2}}%
\sum_{j=2}^{(p-1)/2}\left\{ \frac{1}{j^{2}}+\frac{2\pi ^{2}}{3p^{2}}+\frac{%
\pi ^{4}j^{2}}{9p^{4}}\right\} 
\]%
\[
\leq \frac{p^{2}}{\pi ^{2}}\left\{ \frac{\pi ^{2}}{6}-1+\frac{\pi ^{2}}{3p}+%
\frac{\pi ^{4}}{216p}\right\} <\frac{p^{2}}{\pi ^{2}} 
\]%
if \ $p>10.$ For smaller \ $p$ \ it is easily checked that \ $S(1,p)\geq 0.$

\textbf{(c). }$r\geq 2$ \ We have the expansion%
\[
\cot (x)=\frac{1}{x}-\frac{2}{\pi }\sum_{j=1}^{\infty }\zeta (2j)(\frac{x}{%
\pi })^{2j-1} 
\]%
and hence%
\[
\cot ^{(r)}(x)=\frac{(-1)^{r}r!}{x^{r+1}}\left\{ 1-2(-1)^{r}\sum_{2j\geq
(r+1)}\zeta (2j)\binom{2j-1}{r}(\frac{x}{\pi })^{2j}\right\} 
\]%
We also have the sum%
\[
\sum_{2j\geq r+1}\binom{2j-1}{r}x^{2j}=\frac{x^{r+1}}{2}\left\{ \frac{1}{%
(1-x)^{r+1}}-(-1)^{r}\frac{1}{(1+x)^{r+1}}\right\} 
\]%
and hence%
\[
\left\vert \cot ^{(r)}(\pi /p)\right\vert \geq \frac{r!p^{r+1}}{\pi ^{r+1}}%
\left\{ 1-2\zeta (r+1)\sum_{2j\geq r+1}\binom{2j-1}{r}p^{-2j}\right\} 
\]%
\[
\geq \frac{r!p^{r+1}}{\pi ^{r+1}}\left\{ 1-\zeta (r+1)(\frac{1}{(p-1)^{r+1}}+%
\frac{1}{(p+1)^{r+1}})\right\} 
\]%
In the same way we obtain%
\[
\left\vert \cot ^{(r)}(2\pi /p)\right\vert \leq \frac{r!p^{r+1}}{(2\pi
)^{r+1}}\left\{ 1+\zeta (r+1)((\frac{2}{p-2})^{r+1}+(\frac{2}{p+2}%
)^{r+1})\right\} 
\]%
Furthermore%
\[
\sum_{j=3}^{(p-1)/2}\left\vert \cot ^{(r)}(\pi j/p)\right\vert \leq
\left\vert \dint\limits_{2\pi /p}^{(p-1)\pi /p}\cot ^{(r)}(x)dx\right\vert
\leq \left\vert \cot ^{(r-1)}(2\pi /p)\right\vert 
\]%
\[
\leq \frac{(r-1)!p^{r}}{(2\pi )^{r}}\left\{ 1+\zeta (r)((\frac{2}{p-2})^{r}+(%
\frac{2}{p+2})^{r})\right\} 
\]%
and finally%
\[
\sum_{j=1}^{(p-1)/2}\cot ^{(r)}(\pi j/p)\geq \frac{r!p^{r+1}}{(2\pi )^{r+1}}%
\left\{ 
\begin{array}{c}
1-\zeta (r+1)(\frac{1}{(p-1)^{r+1}}+\frac{1}{(p+1)^{r+1}})-\frac{1}{2^{r+1}}
\\ 
-\zeta (r+1)(\frac{1}{(p-2)^{r+1}}+\frac{1}{(p+2)^{r+1}})-\frac{1}{r2^{r}}
\\ 
-\frac{\zeta (r)}{r}(\frac{1}{(p-2)^{r}}+\frac{1}{(p+2)^{r}})%
\end{array}%
\right\} >0 
\]%
if \ $p\geq 5$ \ (since \ $r\geq 2$ ).

\textbf{Remark 3.1: }Making computer experiments it looks like

\textbf{(a) }\ $T(r,p)$ \ is divisible by \ $p$ \ if not%
\[
r\equiv \frac{p-3}{2}\text{ }\func{mod}\text{ }(p-1) 
\]%
( $c_{p}=T((p-3)/2)$ is not divisible by \ $p$ )

\textbf{(b).} The smallest \ $r$ \ for which \ $T(r,p)$ \ is not an integer
is \ $r=\binom{p}{2}-1.$

\textbf{End of the proof of Proposition 2.2}

In the case \ $p\equiv 3$ $\func{mod}$ $4$ we have \ $\QOVERD( )
{p-j}{p}=-\QOVERD( ) {j}{p}$ \ and formulas (v) and (vi) follow from the
definition of \ $\QOVERD( ) {j}{p}.$

\textbf{Remark 3.2. }Noting that for \ $x\notin \mathbf{Z}$ \ we have \ $%
B_{1}(frac(x))=((x))$ \ we define%
\[
G(p)=T(0,p)=\sum_{j=1}^{p-1}((j^{2}/p)) 
\]%
Then \ $G(p)$ \ is an integer, nonzero if \ $p\equiv 3$ $\func{mod}$ $4$ .
Then computations suggest that \ $G(p)=-1$ if and only if 
\[
p=7,11,19,43,67,163 
\]%
These primes are recognized as those with class number one for \ $\mathbf{Q}(%
\sqrt{-p})$. Indeed, \ $-G(p)$ \ is nothing but Dirichlet's formula for the
class number%
\[
h(-p)=-\frac{1}{p}\sum_{j=1}^{p-1}j\QOVERD( ) {j}{p} 
\]%
which, since \ $p\equiv 3$ $\func{mod}$ $4$ , can be written%
\[
h(-p)=-\sum_{j=1}^{p-1}((j^{2}/p)) 
\]%
or after use of Parseval's formula%
\[
h(-p)=\frac{1}{\sqrt{p}}\sum_{j=1}^{(p-1)/2}\cot (\pi j^{2}/p) 
\]

\textbf{Proof of Conjecture 2.}

We will actually prove a little more:

\textbf{Theorem 3.3. }Let%
\[
f(x)=\dprod\limits_{n=1}^{\infty }\frac{(1-x^{pn})^{p}}{1-x^{n}} 
\]%
Then for \ $(h,k)=(p,k)=1$ \ we have \ ($p\geq 5$ \ prime)%
\[
f(\exp (2\pi ih/k-t) 
\]%
\[
=(\frac{2\pi }{k})^{(p-1)/2}\QOVERD( ) {k}{p}p^{-p/2}t^{(p-1)/2}\exp \left\{ 
\frac{p^{2}-1}{24}t-\frac{\pi i(p^{2}-1)h}{12k}\right\} H(\exp \left\{ \frac{%
2\pi iB}{k}-\frac{4\pi ^{2}}{k^{2}pt}\right\} ) 
\]%
where%
\[
Bph\equiv -1\text{ }\func{mod}\text{ }k 
\]%
and%
\[
H(x)=\frac{\eta (x)^{p}}{\eta (x^{p})}=\dprod\limits_{n=1}^{\infty }\frac{%
(1-x^{n})^{p}}{1-x^{pn}} 
\]

\textbf{Proof: }We have%
\[
f(x)=x^{-(p^{2}-1)/24}\frac{\eta (x^{p})^{p}}{\eta (x)}%
=x^{-(p^{2}-1)/24}h(x) 
\]%
where%
\[
h(x)=\frac{\eta (p\tau )^{p}}{\eta (\tau )}\text{ \ if \ }x=\exp (2\pi i\tau
) 
\]%
Here 
\[
\eta (x)=x^{1/24}\dprod\limits_{n=1}^{\infty }(1-x^{n}) 
\]%
is the Dedekind function. Put \ $\tau =-\frac{1}{p\tau ^{\prime }}$ . Then%
\[
h(\tau )=\frac{\eta (-1/\tau ^{\prime })^{p}}{\eta (-1/p\tau ^{\prime })}=%
\frac{1}{\sqrt{p}}(\frac{\tau ^{\prime }}{i})^{(p-1)/2}H(\tau ^{\prime }) 
\]%
where%
\[
H(\tau ^{\prime })=\frac{\eta (\tau ^{\prime })^{p}}{\eta (p\tau ^{\prime })}
\]%
But by Ogg [8], p.28 \ $H$ \ is a modular form of weight \ $(p-1)/2$ \ for
the group%
\[
\Gamma _{0}(p)=\left\{ \left( 
\begin{array}{cc}
a & b \\ 
c & d%
\end{array}%
\right) \epsilon SL(2,\mathbf{Z);}\text{ }c\equiv 0\text{ }\func{mod}\text{ }%
p\right\} 
\]%
with character \ $\QOVERD( ) {d}{p}$ . This means that%
\[
H(\frac{a\tau ^{\prime }+b}{c\tau ^{\prime }+d})=\QOVERD( ) {d}{p}(c\tau
^{\prime }+d)^{(p-1)/2}H(\tau ^{\prime }) 
\]%
if \ $c\equiv 0$ $\func{mod}$ $p.$ Choose%
\[
\left( 
\begin{array}{cc}
a & b \\ 
c & d%
\end{array}%
\right) =\left( 
\begin{array}{cc}
\frac{1+Bph}{k} & B \\ 
ph & k%
\end{array}%
\right) \epsilon \Gamma _{0}(p) 
\]%
where%
\[
Bhp\equiv -1\text{ }\func{mod}\text{ }k 
\]%
If we take 
\[
\tau =\frac{h}{k}+\frac{2\pi i}{pk^{2}t} 
\]%
and%
\[
\frac{\tau ^{\prime }}{i(c\tau ^{\prime }+d)}=\frac{2\pi }{pkt} 
\]%
Hence%
\[
h(\tau )=\frac{1}{\sqrt{p}}(\frac{2\pi }{pkt})^{(p-1)/2}\QOVERD( )
{d}{p}H(\exp \left\{ \frac{2\pi iB}{k}-\frac{4\pi ^{2}}{k^{2}pt}\right\} ) 
\]%
and the Theorem is proved

Conjecture 2 follows by comparing the limit of this formula and (*) when \ $%
t\rightarrow 0+$%
\[
\]

\textbf{Remark. }$\ $Attempts to find a better approximation to \ $a_{p}(n)$
\ like what is done in [1] failed. The convergence is too slow.

\textbf{Appendix: \ Some Finite Fourier Transformations.}

\textbf{1. }$f(x)=B_{r}(x)$

We have the generating function%
\[
\frac{te^{xt}}{e^{t}-1}=\sum_{r=0}^{\infty }B_{r}(x)\frac{t^{r}}{r!} 
\]%
Putting \ $x=j/k$ \ and summing over \ $j$ \ we get%
\[
\frac{t}{e^{t}-1}\sum_{j=0}^{k-1}\exp (\frac{jt}{k}-\frac{2\pi ij\mu }{k}%
)=\sum_{r=0}^{\infty }\widehat{B}_{r}(\mu /k)\frac{t^{r}}{r!} 
\]%
Now%
\[
\frac{1}{e^{i\alpha }-1}=-\frac{1}{2}-\frac{i}{2}\cot (\alpha /2) 
\]%
so the left hand side is 
\[
\frac{t}{\exp (t/k-2\pi i\mu /k)-1}=t\left\{ -\frac{1}{2}+\frac{i}{2}\cot
(\pi \mu /k+it/2k\right\} 
\]%
\[
=-\frac{1}{2}+k\sum_{r=1}^{\infty }\cot ^{(r-1)}(\pi \mu /k)\frac{t^{r}}{%
(r-1)!}(\frac{i}{2k})^{r}/ 
\]%
Comparing the coefficients of \ $t^{r}$ \ we find%
\[
\widehat{B}_{r}(\mu /k)=kr(\frac{i}{2k})^{r}\cot ^{(r-1)}(\pi \mu /k) 
\]

\textbf{2. \ }$f(x)=l(s,x)$

We have 
\[
\widehat{l}(s,\mu /k)=\sum_{j=0}^{k-1}l(s,j/k)\exp (-2\pi ij\mu /k) 
\]%
\[
=\sum_{n=1}^{\infty }\frac{1}{n^{s}}\sum_{j=0}^{k-1}\exp (2\pi ij(n-\mu
)/k)=k\sum_{n\equiv \mu \text{ }\func{mod}\text{ }k}\frac{1}{n^{s}} 
\]%
\[
=k^{1-s}\sum_{m=0}^{\infty }\frac{1}{(m+\mu /k)^{s}}=k^{1-s}\varsigma (s,\mu
/k) 
\]%
Since 
\[
\widehat{\widehat{f}}(x)=kf(1-x) 
\]%
we have%
\[
\widehat{\zeta }(s,\mu /k)=k^{s}l(s,1-\mu /k) 
\]%
\[
\]

\textbf{References:}

1. G.Almkvist, Asymptotic formulas and generalized Dedekind sums,

\ \  Exp. mathematics, 7, (1998), 343-359

2. G.Almkvist and A.Meurman, Values of Bernoulli polynomials and Hurwitz's 

\ \ \ zeta function at rational points, C.R.Math.Rep.Acad. Canada 13,
(1991), 104-108

3. G.Almkvist, Proof of a conjecture about unimodal polynomials,

\ \  J.Number Theory, 32 (1989), 43-57.

4. T.Apostol, Introduction to analytic number theory, Springer-Verlag,
Berlin 1976.

5. T.Apostol, Modular functions and Dirichlet series in number theory,

\ \  Springer-Verlag, Berlin 1976.

6. Z.I.Borevich and I.R.Shafarevich, Number theory, Acad. Press, New York
1966.

7. T. Dokshitzer, On Wilf's conjecture and generalizations, pp.133-153 in
Number Theory 

\ \ \ (Halifax, Canada, I,1994),  edited by K.Dilcher, CMS Conf. Proc. 15,

\ \  Amer. Math. Soc. Providence, RI, 1995

8. F.Garvan, Some congruences for partitions that are p-cores,

\ \  Proc. London Math. Soc. 66 (1993), 449-478.

9. A.Ogg, Survey of modular functions of one variable,

\ \  Springer Lecture Notes, Vol 320, p.1-36.

10. S.Ramanujan, On certain trigonometric sums and their applications 

\ \ \ \ in the theory of numbers, Trans.Cambridge Phil. Soc. 22 (1918),
259-276.

\bigskip 

Institute of Algebraic Meditation

Fogdar\"{o}d 208

S-24333 H\"{o}\"{o}r, Sweden

gert@maths.lth.se

\end{document}